\documentclass{amsart}
\usepackage{amsmath,amssymb}
\newtheorem{theorem}{Theorem}[section]
\newtheorem{lemma}[theorem]{Lemma}
\newtheorem{proposition}[theorem]{Proposition}
\newtheorem{corollary}[theorem]{Corollary}

\newtheorem{remark}[theorem]{Remark}
\newtheorem{example}[theorem]{Example}

\begin{document}
\title{Indecomposable modules and Gelfand rings}
\author{Fran\c cois Couchot}
\address{Laboratoire de Math\'ematiques Nicolas Oresme, CNRS UMR
  6139,
D\'epartement de math\'ematiques et m\'ecanique,
14032 Caen cedex, France}
\email{couchot@math.unicaen.fr} 


\begin{abstract}  It is proved that a commutative ring is clean if and only if it is Gelfand
with a totally disconnected maximal spectrum. It is shown that each indecomposable module over a commutative ring $R$ satisfies a finite condition if and 
only if $R_P$ is an artinian valuation ring for each maximal prime ideal $P$. Commutative rings for which each
indecomposable module has a local endomorphism ring are studied. These rings
are clean and elementary divisor rings. It is shown that each commutative ring $R$ with a Hausdorff
and totally disconnected maximal spectrum is local-global.
Moreover, if $R$ is arithmetic then $R$ is an elementary divisor
ring. 
\end{abstract}

\maketitle

In this paper $R$ is a  commutative ring with unity and modules are
unitary.

In \cite[Proposition 2]{GoWa76} Goodearl and Warfield proved that each
zero-dimensional ring $R$ satisfies the second condition of our
Theorem~\ref{T:geldis}, and this condition plays a crucial role in
their paper. In Section~\ref{S:logo}, we  show that a ring $R$ enjoys
this condition if and only if it is clean, if and only
if it is Gelfand with a totally disconnected maximal spectrum. So we get a
generalization of results obtained by Anderson and Camillo in
\cite{AnCa02} and by Samei in \cite{Sam04}. We deduce that every
commutative ring $R$ with a Hausdorff and totally disconnected maximum
prime spectrum is local-global, and moreover, $R$ is an elementary divisor ring
if, in addition, $R$ is arithmetic.
One can see in \cite{EsGu82} that local-global rings have very
interesting properties.

In Section~\ref{S:maxi} we give a characterization of commutative rings for which each indecomposable module satisfies a finite condition: finitely 
generated, finitely cogenerated, cyclic, cocyclic, artinian, noetherian or of finite length. We deduce  that a commutative ring is 
Von Neumann regular if and only if each indecomposable module is
simple. This last result was already proved in \cite{Couc03}. We study
commutative rings for which each indecomposable module has a local
endomorphism ring. These rings are clean and elementary divisor
rings. It remains to find valuation rings satisfying this property to give a complete characterization of
these rings. We also give characterizations 
of  Gelfand rings and clean rings by using properties of indecomposable modules. Similar results are obtained in Section~\ref{S:mini}, for commutative rings for which
each prime ideal contains only one minimal prime ideal.

We denote respectively $\mathrm{Spec}\ R$, $\mathrm{Max}\ R$ and $\mathrm{Min}\ R,$ the
space of prime ideals, maximal ideals, and minimal prime ideals of
$R$, with the Zariski topology. If $A$ a subset of $R$, then we denote  
\[V(A) = \{ P\in\mathrm{Spec}\ R\mid A\subseteq P\},\ D(A) = \{ P\in\mathrm{Spec}\ R \mid A\not\subseteq P\},\]
\[V_M(A)=V(A)\cap\mathrm{Max}\ R\ \mathrm{and}\ D_M(A) =D(A)\cap\mathrm{Max}\ R.\]

\section{Local-global Gelfand rings} \label{S:logo}

As in \cite{Mul79} we say that a commutative ring $R$ is \textbf{Gelfand} if each prime ideal is contained
in only one maximal ideal. In this case,  we put $\mu:\mathrm{Spec}\ R\rightarrow\mathrm{Max}\ R$
the map defined by $\mu(J)$ is the maximal ideal containing $J$ for each prime ideal $J$. Then
$\mu$ is continuous and $\mathrm{Max}\ R$ is Hausdorff by \cite[Theorem 1.2]{MaOr71}.

In \cite{GoWa76}, Goodearl and Warfield proved that every zero-Krull-dimensional commutative ring satisfies the second condition of the following 
theorem. This property is used to show cancellation, $n$-root and isomorphic refinement theorems for finitely generated modules
over algebras over a commutative ring which is Von Neumann regular modulo its Jacobson radical. So, the following theorem allows
us to extend these results to each ring with a Hausdorff and totally disconnected
maximal spectrum. As in \cite{Nic77} we say that a ring $R$ is \textbf{clean} if each
element of $R$ is the sum of a unit with an idempotent. In
\cite[Proposition 1.8 and Theorem 2.1]{Nic77} Nicholson proved that
commutative clean rings are exactly the exchange rings defined by
Warfield in \cite{War72}. In \cite{Sam04}
Samei proved that the conditions (1), (3) and (4) are equivalent when $R$
is semiprimitive and in \cite{AnCa02} Anderson and Camillo showed that each clean ring
is Gelfand. We can also see \cite[Theorem 3]{Mon72}. If $P$ is a prime ideal we denote by $0_P$ the kernel of
the natural map $R\rightarrow R_P$.

\begin{theorem} \label{T:geldis} Let $R$ be a ring. The
  following conditions are equivalent:
\begin{enumerate} 
\item $R$ is a Gelfand ring and $\mathrm{Max}\ R$ is totally disconnected.
\item Each $R$-algebra $S$ satisfies this condition: 
let
  $f_1,\dots,f_k$ be polynomials over $S$ in noncommuting variables
  $x_1,\dots,x_m,y_1,\dots,y_n$. Let $a_1,\dots,a_m\in S$. Assume that
   $\forall P\in \mathrm{Max}\ R$ there exists $b_1,\dots,b_n\in S_P$
  such that $f_i(a_1,\dots,a_m,b_1,\dots,b_n)=0$  $\forall i$, $1\leq
  i\leq k$. Then there exist $d_1,\dots,d_n\in S$ such that
  $f_i(a_1,\dots,a_m,d_1,\dots,d_n)=0$  $\forall i$, $1\leq
  i\leq k$. 
\item $R$ is a clean ring.
\item $R$ is Gelfand and $\forall P\in\mathrm{Max}\ R,\ 0_P$ is
  generated by  a set of idempotents.
\end{enumerate}

\end{theorem}
\textbf{Proof.} $(1)\Rightarrow (2)$. By \cite[Theorem 16.17]{GiJe60}
$\mathrm{Max}\ R$ has a base of clopen subsets. Since $\mu$ is continuous, each
clopen subset of $\mathrm{Max}\ R$ is of the form $D_M(e)$ for some
idempotent $e$. So we can do the same proof as in
\cite[Proposition 2]{GoWa76} where we replace $\mathrm{Spec}\ R$ with $\mathrm{Max}\ R$.

$(2)\Rightarrow (3)$. Let $a\in R$. We consider the following
equations: $x^2=x$ and $y(a-x)=1$. Since each local ring is clean,
these equations have a solution in $R_P$ for each maximal ideal
$P$. We conclude that there is also a solution in $R$ and that $R$ is clean.

$(3)\Rightarrow (1)$. Let $P,P'\in\mathrm{Max}\ R$, $P\ne P'$. Then there
exist $a\in P$ and $a'\in P'$ such that $a+a'=1$. We have $a=u+e$
where $u$ is a unit and $e$ is an idempotent. Since $a\in P$ and
$u\notin P$ we get that $e\notin P$. We have $a'=1-a=-u+(1-e)$. So
$1-e\notin P'$. Consequently $P$ and $P'$ have disjoint clopen
neighbourhoods. Since $\mathrm{Max}\ R$ is
quasi-compact, we deduce that this space is compact and totally
disconnected. The equality $e(1-e)=0$ implies that $P\cap P'$ contains
no prime ideal. Hence $R$ is Gelfand.

$(1)\Rightarrow (4)$. Let $P$ be a maximal ideal and $a\in 0_P$. Then
there exists $s\in R\setminus P$ such that $sa=0$. Since
$\mathrm{Max}\ R$ is totally disconnected there is a clopen subset $U$
such that $U\subseteq D_M(s)$. Because of $\mu$ is continuous,
there exists an idempotent $e$ such that $P\in
D(e)=\mu^{\leftarrow}(U)\subseteq\mu^{\leftarrow}(D_M(s))\subseteq
D(s).$ Then $e\in Rs$, $ea=0$, $a=a(1-e)$ and $1-e\in 0_P$.

$(4)\Rightarrow (1)$. Let $P,P'\in\mathrm{Max}\ R$, $P\ne P'$. Since
 $R$ is Gelfand, by \cite[Theorem 1.2]{MaOr71} there exist $a\in
 0_P\setminus P'$. Then there exists an
 idempotent $e\in 0_P\setminus P'$. Clearly $1-e\notin P$.
 Consequently $P$ and $P'$
have disjoint clopen neighbourhoods.
 \qed

\bigskip
We say that $R$ is \textbf{local-global} if each polynomial over $R$ in
finitely many indeterminates which admits unit values locally, admits
unit values. Recall that most of the results of \cite{GoWa76} about commutative rings which are Von Neumann regular modulo 
their Jacobson radicals, have been extended to local-global rings by
Estes and Guralnick in \cite{EsGu82}. We have the following corollary:
\begin{corollary} \label{C:loglo} Let $R$ be a ring such that
  $\mathrm{Max}\ R$ is Hausdorff and totally disconnected. Then $R$ is
  local-global. 
\end{corollary}
\textbf{Proof.} Let $J$ be the Jacobson radical of $R$. Then $R$ is
local-global if and only if $R/J$ is local-global. So we may assume
that $R$ is semiprimitive. From the remark that follows \cite[Theorem
1.2]{MaOr71} and from Theorem~\ref{T:geldis} we deduce that $R$ is clean. Let $f$ be a polynomial over $R$ in
finitely many indeterminates $X_1,\dots,X_n$, which admits unit values
locally. Then, we apply theorem~\ref{T:geldis} by taking $S=R$ to the
polynomial $Yf(X_1,\dots,X_n)-1$.  \qed

\begin{remark} \textnormal{If $R$ is the ring of algebraic integers, then
$R$ is local-global by \cite{Dad64} and semi-primitive. But this ring
is not Gelfand.} 
\end{remark}

\bigskip

\section{Arithmetic Gelfand rings} \label{S:ari}

We say that a module is \textbf{uniserial} if its set of submodules is totally ordered
by inclusion, we say
that a ring $R$ is a \textbf{valuation ring} if it is uniserial as $R$-module and we say that $R$ is \textbf{arithmetic} if $R_P$ is a valuation ring for
each maximal ideal $P$.
Recall that $R$ is a \textbf{B\'ezout ring} if each finitely generated ideal is principal and 
$R$ is an \textbf{elementary divisor ring} if each finitely presented module is a direct sum of cyclic submodules.
\begin{theorem} \label{T:ariloglo} Let $R$ be an arithmetic local-global ring
 . Then $R$ is an elementary divisor ring. Moreover, for each $a,b\in R$, there
  exist $d,a',b',c\in R$ such that $a=a'd$, $b=b'd$ and 
  $a'+cb'$ is a unit of $R$. 
\end{theorem}
\textbf{Proof.}  Since every finitely generated ideal is locally
principal $R$ is B\'ezout by \cite[Corollary 2.7]{EsGu82}. Let $a,b\in
R$. Then there exist $a',b',d\in R$ such that $a=a'd,b=b'd$ and
$Ra+Rb=Rd$. Consider the following polynomial $a'+b'T$. If $P$ is a
maximal ideal, then we have $aR_P=dR_P$ or $bR_P=dR_P$. So, either $a'$ is a
unit of $R_P$ and $a'+b'r$ is a unit of $R_P$ for each $r\in PR_P$, or
$b'$ is a unit of $R_P$ and $a'+b'(1-a'/b')$ is a unit of $R_P$. We
conclude that the last assertion holds. Now, let $a,b,c\in R$ such that
$Ra+Rb+Rc=R$. We set $Rb+Rc=Rd$. Let $b',c',s$ and $q$ such that
$b=b'd,\ c=c'd$ and $b'+c'q$ and $a+sd$ are units. Then
$(b'+c'q)(a+sd)=(b'+c'q)a+s(b+qc)$ is a unit. We conclude by
\cite[Theorem 6]{GiHe56}.  \qed

\bigskip
We deduce the following corollary which is a generalization of \cite[Theorem
III.6]{Cou03} and \cite[Theorem 5.5]{GiHen56}.
\begin{corollary} \label{C:arigel} Let $R$ be an arithmetic  ring
  with a  Hausdorff and totally disconnected maximal spectrum. Then $R$ is an
  elementary divisor ring. Moreover, for each $a,b\in R$, there
  exist $d,a',b',c\in R$ such that $a=a'd$, $b=b'd$ and
  $a'+cb'$ is a unit of $R$. 
\end{corollary}

\begin{corollary} Let $R$ be an arithmetic Gelfand ring
  such that $\mathrm{Min}\ R$ is compact. Then $R$ is an
  elementary divisor ring. Moreover, for each $a,b\in R$, there
 exist $d,a',b',c\in R$ such that $a=a'd$, $b=b'd$ and
  $a'+cb'$ is a unit of $R$.
\end{corollary}
\textbf{Proof} Let $\mu'$ be the restriction of $\mu$ to $\mathrm{Min}\ R$.
Since $R$ is arithmetic each prime ideal contains only one minimal prime ideal.  
Then $\mu'$ is bijective and it is an homeomorphism because $\mathrm{Min}\ R$ is compact.
One can apply
corollary~\ref{C:arigel} since
$\mathrm{Min}\ R$ is totally disconnected . \qed

\begin{remark} \textnormal{In \cite{GiHen56} there is an example of a Gelfand B\'ezout ring $R$ which is not
an elementary divisor ring. Consequently $\mathrm{Min}\ R$
is not compact.}
\end{remark}

\section{Indecomposable modules and maximal ideals} \label{S:maxi}

In the two next propositions we give a characterization of Gelfand
rings and clean rings by using properties of indecomposable modules.

\begin{proposition} \label{P:Gsupp}
Let $R$ be a ring. The following conditions are
  equivalent:
\begin{enumerate}
\item For each $R$-algebra $S$ and for each left $S$-module $M$ for which $\mathrm{End}_S(M)$ is local,
  $\mathrm{Supp}\ M$ contains only one maximal ideal.
\item $R$ is a Gelfand ring.
\item $\forall P\in\mathrm{Max}\ R$ the natural map $R\rightarrow R_P$ is surjective.
\end{enumerate}

When these conditions are satisfied, $M=M_P$ for each left $S$-module $M$ for which
$\mathrm{End}_S(M)$ is local, where $P$ is the unique
maximal ideal of $\mathrm{Supp}\ M$ and where $S$ is an algebra over $R$.
\end{proposition}
\emph{Proof.} Assume that $R$ is Gelfand. Let $S$ be an $R$-algebra
 and let $M$ be a left $S$-module such that $\mathrm{End}_S(M)$ is local.
 Let $P$ be the prime ideal which is the inverse image of
the maximal ideal of $\mathrm{End}_S(M)$ by the canonical map
$R\rightarrow\mathrm{End}_S(M)$ and let $Q=\mu(P)$.
 Since $M$ is an $R_P$-module,
$0_Q\subseteq\mathrm{ann}_R(M)$. So, $\mathrm{Supp}\ M\subseteq V(0_Q)$ and $Q$ is the only maximal ideal
belonging to $V(0_Q)$ since $R$ is Gelfand.

Conversely, if $P$ is a prime ideal then $R_P=\mathrm{End}_R(R_P)$. It
follows that $P$ is contained in only one maximal ideal. 

By \cite[Theorem 1.2]{MaOr71} $R$ is Gelfand if and only if, $\forall P\in\mathrm{Max}\ R$, $P$ is the 
only maximal ideal containing $0_P$. This is equivalent to $R/0_P$ is local, $\forall P\in\mathrm{Max}\ R$.
It is obvious that $R_P=R/0_P$ if $R/0_P$ is local. (When $R$ is semi-primitive we can see
\cite[Proposition 1.6.1]{Bko70}).

Recall that the diagonal map $M\rightarrow \Pi_{P'\in\mathrm{Max}(R)}M_{P'}$
  is monic. Since $R$ is Gelfand, we have $M_P=M/0_PM$ where $P$ is the only maximal ideal of
$\mathrm{Supp}\ M$. Hence $M=M_P$. 
\qed

\begin{proposition} \label{P:Gind}
Let $R$ be a ring. The following conditions are
  equivalent:
\begin{enumerate}
\item For each $R$-algebra $S$ and for each indecomposable left $S$-module $M$,
$\mathrm{Supp}\ M$ contains only one maximal ideal.
\item $R$ is clean.
\end{enumerate}
When these conditions are satisfied, $M=M_P$ for each indecomposable left
$S$-module $M$, where $P$ is the unique maximal ideal of
$\mathrm{Supp}\ M$ and $S$ is an $R$-algebra.
\end{proposition}
\textbf{Proof.} $(2)\Rightarrow (1).$  By Theorem~\ref{T:geldis}
$\mathrm{Max}\ R$ is totally disconnected. So, if $P$ and $P'$ are two different maximal ideals
such that $P\in\mathrm{Supp}\ M$ then there exists an idempotent $e\in
P\setminus P'$ because $\mu$ is continuous. Since $(1-e)\notin P$ and $M_P\not=0$, we have
$(1-e)M\not=0$. We deduce that $eM=0$ and $M_{P'}=0$.

$(1)\Rightarrow (2).$ $R$ is Gelfand by Proposition~\ref{P:Gsupp}. Let $A$ be an
ideal such that $V(A)$ is the inverse image of a connected component of
$\mathrm{Max}\ R$ by $\mu$. Then $V(A)$ is connected too, whence $R/A$ is indecomposable. So there is only one
maximal ideal in $V(A)$. Since each
connected component contains only one point, 
$\mathrm{Max}\ R$ is totally disconnected. \qed

\bigskip

This lemma is needed to prove the main results of this section.
\begin{lemma} \label{L:val} Let $R$ be a local ring which is not a
  valuation ring. Then there exists an indecomposable  non-finitely
  generated $R$-module whose endomorphism ring is not local.
 \end{lemma}
\textbf{Proof.}
Since $R$ is not a valuation ring there exist $a,\ b\in R$ such that neither divides
the other. By taking a suitable quotient ring, we may assume that
$Ra\cap Rb=0$ and $Pa=Pb=0$. Let $F$ be a free module generated by
$\{e_n\mid n\in\mathbb{N}\}$, let $K$ be the submodule generated by
$\{ae_n-be_{n+1}\mid n\in\mathbb{N}\}$ and let $M=F/K$. Clearly
$M/PM\cong F/PF$. We will show that $M$ is indecomposable and
$S:=\mathrm{End}_R(M)$ is not local.
Let us observe that $M$ is defined as in proof of \cite[Theorem
2.3]{Gri70}. But, since $R$ is not necessarily artinian, we do a
different proof to show that $M$ is indecomposable. We shall prove
that $S$ contains no trivial idempotents. Let $s\in S$. Then $s$ is
induced by an endomorphism $\tilde{s}$ of $F$ which satisfies
$\tilde{s}(K)\subseteq K$. For each $n\in\mathbb{N}$ there exists a
finite family $(\alpha_{p,n})$ of elements of $R$ such that:
\begin{equation} \label{eq:s1}
\tilde{s}(e_n)=\sum_{p\in\mathbb{N}}\alpha_{p,n}e_p
\end{equation}
Since $\tilde{s}(K)\subseteq K,\ \forall n\in\mathbb{N},\ \exists$ a
finite family $(\beta_{p,n})$ of elements of $R$ such that:
\begin{equation} \label{eq:s2}
\tilde{s}(ae_n-be_{n+1})=\sum_{p\in\mathbb{N}}\beta_{p,n}(ae_p-be_{p+1})
\end{equation}
From \ref{eq:s1} and \ref{eq:s2} 
if follows that:
\[\sum_{p\in\mathbb{N}}(a\alpha_{p,n}-b\alpha_{p,n+1})e_p=a\beta_{0,n}e_0+\sum_{p\in\mathbb{N}^*}(a\beta_{p,n}-b\beta_{p-1,n})e_p\]
Since $Pa=Pb=Ra\cap Rb=0$ we deduce that 
\[\alpha_{0,n+1}\equiv 0\ [P],\ \ \alpha_{p,n}\equiv\beta_{p,n}\ [P]\ \ 
{\rm and}\ \ 
\alpha_{p,n+1}\equiv\beta_{p-1,n}\ [P]\]
It follows that 
\begin{equation} \label{eq:s3}
(i)\ \alpha_{p,n}\equiv\alpha_{p+1,n+1}\ [P],\ \forall
p,n\in\mathbb{N},\ \
{\rm and}\ \ 
(ii)\ \alpha_{p,p+k+1}\equiv 0\ [P],\ \forall p,k\in\mathbb{N}
\end{equation}

Now we assume that $s$ is idempotent. Let $x_n=e_n+K,\ \forall
n\in\mathbb{N}$. Let $\bar{s}$ be the
endomorphism of $M/PM$ induced by $s$. If $L$ is an $R$-module and $x$
an element of $L$, we put $\bar{x}=x+PL$. From $s^2(x_0)=s(x_o)$ we
get the following equality:
\begin{equation} \label{eq:s5}
\sum_{n\in\mathbb{N}}(\sum_{p\in\mathbb{N}}\alpha_{n,p}\alpha_{p,0})x_n=\sum_{n\in\mathbb{N}}\alpha_{n,0}x_n
\end{equation}
Then
$\bar{\alpha}_{0,0}=\sum_{p\in\mathbb{N}}\bar{\alpha}_{0,p}\bar{\alpha}_{p,0}=\bar{\alpha}_{0,0}^2$,
since $\bar{\alpha}_{0,p}=0$ by \ref{eq:s3}$(ii)$, $\forall p>0$. So, we
have $\bar{\alpha}_{0,0}=0$ or $\bar{\alpha}_{0,0}=1$. If
$\bar{\alpha}_{0,0}=1$ then we replace $s$ with $\mathbf{1}_M-s$. So
we may assume that $\bar{\alpha}_{0,0}=0$. By \ref{eq:s3}$(i)$ 
 $\bar{\alpha}_{n,n}=0$, $\forall n\in\mathbb{N}$. By using
\ref{eq:s5} and \ref{eq:s3}$(ii)$ we get that 
\[ \bar{\alpha}_{n,0}=\sum_{p=0}^{n-1}\bar{\alpha}_{n,p}\bar{\alpha}_{p,0}+\bar{\alpha}_{n,n}\bar{\alpha}_{n,0}\]
Hence, if $\bar{\alpha}_{p,0}=0,\ \forall p<n$ then
$\bar{\alpha}_{n,0}=0$ too. By induction we obtain that
$\bar{\alpha}_{n,0}=0$,
$\forall n\in\mathbb{N}$. We deduce that
\begin{equation} \label{eq:s6}
\alpha_{p,n}\in P,\ \forall p,n\in\mathbb{N}
\end{equation}

Let $A=\mathrm{Im}\ s,\ B=\mathrm{Ker}\ s$ and let $A'$ and $B'$ be the
inverse image of $A$ and $B$ by the natural map $F\rightarrow M$. If
$x\in A'$ then $\tilde{s}(x)=x+y$ for some $y\in K$. By \ref{eq:s6}
and $Pa=Pb=0$ it follows that $\tilde{s}(y)=0$ and
$\tilde{s}^2(x)=\tilde{s}(x)$. If $x\in B'$ then $\tilde{s}(x)\in
K$. So $\tilde{s}^2(x)=0$. We deduce that
$(\tilde{s}^2)^2=\tilde{s}^2$. Let $C=\mathrm{Im}\ \tilde{s}^2$. Then
$C$ is projective and $C=PC$ by \ref{eq:s6}. By
\cite[Proposition 2.7]{Bas60} $C=0$ . So $s=0$ (or $\mathbf{1}_M-s=0$).

It remains to prove that $S$ is not local. Let $f,g\in S$ defined in
the following way: $f(x_n)=x_{n+1}$ and
$g(x_n)=x_n-x_{n+1}$, $\forall n\in\mathbb{N}$. We
easily check that $x_0\notin \mathrm{Im}\ f\cup \mathrm{Im}\ g$.
So $f$ and $g$ are not units of $S$ and $f+g=\mathbf{1}_M$ is a unit.
Hence $S$ is not local.
\qed

\bigskip
 A module is \textbf{cocyclic} (respectively \textbf{finitely cogenerated})
if it is a submodule of the injective hull of a simple module
(respectively of a finite direct sum of
injective hulls of simple modules).

Now we give a characterization of commutative rings for which each indecomposable module satisfies a 
finite condition.
\begin{theorem} \label{T:ind}
Let $R$ be a ring. The following conditions are
  equivalent:
\begin{enumerate}
\item Each indecomposable $R$-module is of finite length.
\item Each indecomposable $R$-module is noetherian.
\item Each indecomposable $R$-module is finitely generated.
\item Each indecomposable $R$-module is artinian.
\item Each indecomposable $R$-module is finitely cogenerated.
\item Each indecomposable $R$-module is cyclic.
\item Each indecomposable $R$-module is cocyclic.
\item For each maximal ideal $P$, $R_P$ is an artinian valuation ring.
\item $R$ is an arithmetic ring of Krull-dimension $0$ and its
  Jacobson ideal $J$ is T-nilpotent. 
\end{enumerate}
\end{theorem}
\textbf{Proof.} The following implications are obvious: $(1)\Rightarrow
(2)\Rightarrow (3)$, $(1)\Rightarrow (4)\Rightarrow (5)$, $(6)\Rightarrow (3)$ and
$(7)\Rightarrow (5)$.

$(8)\Rightarrow (1),\ (6)\ \mathrm{and}\ (7)$. $R$ has Krull dimension
0. Hence $R$ is clean. So, if $M$ is an indecomposable module, by
proposition~\ref{P:Gind} there is only one maximal ideal $P$ such that
$M_P\not=0$. Moreover $M\cong M_P$. We conclude by \cite[Theorem 4.3]{Gri70}.

$(3)\Rightarrow (8)$. Let $P$ be a maximal ideal and $E$ the injective hull of
$R/P$. Then each submodule of $E$ is indecomposable. It follows that
$E$ is a noetherian module. By \cite[Proposition 3]{Vam68} $E$ is a
module of finite length, and by \cite[Theorem 3]{Vam68} $R_P$ is
artinian. We conclude by \cite[Theorem 2.3]{Gri70} or Lemma~\ref{L:val}.

$(5)\Rightarrow (8)$. Let $P$ be a maximal ideal. Then each factor of $R_P$
modulo an ideal of $R_P$ is finitely cogenerated. It follows that
$R_P$ is artinian. We conclude as above.

$(8)\Rightarrow (9)$. Let $x_1,\dots,x_n,\dots$ be a sequence of elements
of $J$. Then for each maximal ideal $P\ \exists s_P\notin P$ and a
positive integer $n_P$ such that $s_Px_1\dots x_{n_P}=0$. There is a
finite family of open sets
$D(s_{n_{P_1}}),\dots,D(s_{n_{P_m}})$ that cover $\mathrm{Spec}\ R$. We
set $n=max\{n_{P_1},\dots,n_{P_m}\}$. Then $x_1\dots x_n=0$.

$(9)\Rightarrow (8)$. $\forall P\in\mathrm{Max}\ R$, $R_P$ is a
valuation ring and $PR_P$ is a nilideal. Then for every
$r\in P$ there exists $s\in R\setminus P$ such that $sr$ is
nilpotent. So we get that $PR_P=JR_P$, whence $PR_P$ is
T-nilpotent. We easily prove that $R_P$ is artinian. \qed

\bigskip
From this theorem it is easy to deduce the two following
corollaries. Another proof of the second corollary is given in \cite[Theorem
2.13]{Couc03}. 
\begin{corollary} Let $n$ be a positive integer, $R$ a ring and $J$ its
  Jacobson radical. Then
  the following conditions are equivalent:
\begin{enumerate}
\item Each indecomposable module has a length $\leq n$.
\item For each maximal ideal $P$, $R_P$ is a valuation ring and
  $(PR_P)^n=0$. 
\item $R$ is an arithmetic ring of Krull-dimension $0$ and $J^n=0$.
\end{enumerate}
\end{corollary}

\begin{corollary} A ring $R$ is Von Neumann regular if and only if
  every indecomposable module is simple.
\end{corollary}

The next theorem gives a partial characterization of commutative rings
for which each indecomposable module has a local endomorphism ring.
\begin{theorem} \label{C:locend}
Let $R$ be a ring for which $\mathrm{End}_R(M)$ is local for each
indecomposable module $M$. Then
$R$ is a clean elementary divisor ring.

\end{theorem}
\textbf{Proof.} Let $P$ be a prime ideal. Then
$R/P=\mathrm{End}_R(R/P)$ is local. Hence $R$ is Gelfand. We prove
that $\mathrm{Max}\ R$ is totally disconnected as in proof of
 proposition~\ref{P:Gind}. If $P$ is a maximal ideal, each
 indecomposable $R_P$-module $M$ is also indecomposable over $R$ and
 $\mathrm{End}_R(M)=\mathrm{End}_{R_P}(M)$. By Lemma~\ref{L:val} $R_P$
 is a valuation ring. \qed

\begin{example} \textnormal{ If $R$ is a ring satisfying the
    equivalent conditions of Theorem~\ref{T:ind}, then each
    indecomposable $R$-module has a local endomorphism ring. But, by
    \cite[Corollary 2 p.52]{Kap69} and \cite[Corollary 3.4]{ShLe74}, each
    complete discrete rank one valuation ring enjoys this property
    too. So, we consider a complete discrete rank one valuation ring
    $D$, $Q$ its ring of fractions and $R$ the subring of
    $Q^{\mathbb{N}}$ defined as in \cite[Example 1.7]{Nic77}:
    $x=(x_n)_{n\in\mathbb{N}}\in R$ if $\exists p\in\mathbb{N}$
    and $s\in D$ such that $x_n=s,\ \forall n>p$. Since $D$ is local,
    $R$ is clean and semi-primitive by \cite[Theorem 2]{Mon72}. We put
    $\mathbf{1}=(\delta_{n,n})_{n\in\mathbb{N}}$ and $\forall
    p\in\mathbb{N},\ \mathbf{e}_p=(\delta_{p,n})_{n\in\mathbb{N}}$
    where $\delta_{n,p}$ is the Kronecker symbol. Let $J$ be the
    maximal ideal of $D$. If $P$ is a maximal ideal of $R$, then either
    $\mathbf{e}_p\in P,\ \forall p\in\mathbb{N},$ whence
   $P=J\mathbf{1}+\oplus_{p\in\mathbb{N}}R\mathbf{e}_p$ and
    $R_P\cong R/\oplus_{p\in\mathbb{N}}R\mathbf{e}_p\cong D$, or
    $\exists p\in\mathbb{N}$ such that $\mathbf{e}_p\notin P,$ whence
   $P=R(\mathbf{1}-\mathbf{e}_p)$ and $R_P\cong R/P\cong Q$. Thus $R$ is arithmetic and each
    indecomposable $R$-module has a local endormorphism ring. Observe
    that each indecomposable $R$-module is uniseriel and linearly
    compact and its endomorphism ring is commutative.}
\end{example}

\section{Indecomposable modules and minimal prime ideals}
\label{S:mini}
In this section we study rings $R$ for which each prime ideal contains
only one minimal prime ideal. 
In this case, if $P\in\mathrm{Spec}\ R$, let $\lambda(P)$ be the only
minimal prime ideal contained 
in $P$. We shall see that $\lambda$ is continuous if and only if
$\mathrm{Min}\ R$ is compact.
(See \cite[Theorem 2]{Kis74} when $R$ is semi-prime). But, since $\lambda$ is surjective, the 
set of minimal primes can be endowed with the quotient topology induced by the Zariski topology
of $\mathrm{Spec}\ R$. We denote this topologic space by $\mathrm{QMin}\ R$. Then we have the
following:

\begin{proposition} \label{P:pSHau} Let $R$ be a ring such that each
  prime ideal contains a unique minimal prime ideal and $N$ its
  nilradical. Then $\mathrm{QMin}\ R$ is compact. Moreover, $\mathrm{QMin}\ R$ and
$\mathrm{Min}\ R$ are homeomorphic if and only if $\mathrm{Min}\ R$ is compact.
\end{proposition}
The following lemma is needed to prove this proposition. This lemma is
a ge\-neralization of \cite[Lemma 2.8]{HeJe65}. We do a  similar proof.
\begin{lemma} \label{L:pri} Let $R$ be a ring, $N$ its nilradical and
  $a\in R\setminus N$. Let $P$ be a prime ideal such that $P/(N:a)$ is minimal in
  $R/(N:a)$. Then $P$ is a minimal prime ideal.
\end{lemma}
\textbf{Proof.} First we show that $a+(N:a)$ is a non-zerodivisor in $R/(N:a)$ and consequently
$a\notin P$. Let $b\in R$
such that $ab\in (N:a)$. Then $a^2b\in N$. We easily deduce that $ab\in N$, whence $b\in (N:a)$.
Let $r\in P$. Then there exist a positive integer $n$ and $s\in R\setminus P$ such that $sr^n\in (N:a)$.
It follows that $asr^n\in N$. Since $as\notin P$ we deduce that $PR_P$ is a nilideal, whence $P$ 
is a minimal prime. \qed

\bigskip
\textbf{Proof of proposition~\ref{P:pSHau}.} Let $A$ and $B$ be two distinct
 minimal prime ideals. Since each maximal
ideal contains only one minimal prime ideal, we have $A+B=R$. Therefore
there exist $a\in A$ and $b\in B$ such that $a+b=1$. Thus $a\notin
B$ and $a\notin N$. But $a$ is a nilpotent element of $R_A$. Hence
$(N:a)\nsubseteq A$. In the same way we show that $B\in
D((N:b))$. We have $(N:a)\cap (N:b)=(N:Ra+Rb)=N$. So $D((N:a))\cap
D((N:b))=\emptyset$. By Lemma~\ref{L:pri}, $D((N:a))$ and $D((N:b))$ are the inverse images of
disjoint open subsets of $\mathrm{QMin}\ R$ by $\lambda$. We conclude that this space
 is Hausdorff. Since $\mathrm{Spec}\ R$ is quasi-compact,
it follows that $\mathrm{QMin}\ R$ is compact. 

Let $\lambda'$ be the restriction of $\lambda$ to $\mathrm{Min}\ R$. It is obvious
that $(\lambda')^{-1}$ is continuous if and only if $\mathrm{Min}\ R$ is compact.  
\qed  

\begin{remark} \textnormal{If we consider the set of D-components of $\mathrm{Spec}\ R$, defined in \cite{Laz67},
endowed with the quotient topology, we get a topologic space $X$. Then $X$ is homeomorphic to $\mathrm{Max}\ R$
(respectively $\mathrm{QMin}\ R$) if $R$ is Gelfand (respectively every prime ideal contains only one minimal
prime). But $X$ is not generally Hausdorff: see \cite[Propositions 6.2 and 6.3]{Laz67}.}
\end{remark}

\bigskip Now we can show the two following propositions which are similar to Propositions~\ref{P:Gsupp}
and \ref{P:Gind}. The proofs are similar too.
\begin{proposition} \label{P:Bsupp}
Let $R$ be a ring. The following conditions are
  equivalent:
\begin{enumerate}
\item For each $R$-algebra $S$ and for each left $S$-module $M$ for which $\mathrm{End}_S(M)$ is local, there exists
 only one minimal prime ideal $A$ such that $\mathrm{Supp}\ M\subseteq V(A).$
\item Every prime ideal contains only one minimal prime ideal.
\end{enumerate}
\end{proposition}
\textbf{Proof.}  $(2)\Rightarrow (1)$. Let $S$ be an $R$-algebra and
let $M$ be a left $S$-module such that $\mathrm{End}_S(M)$ is local.
Let $P$ be the prime ideal which is the inverse image of
the maximal ideal of $\mathrm{End}_S(M)$ by the canonical map
$R\rightarrow\mathrm{End}_S(M)$, $A=\lambda(P)$ and $0_P$ the kernel of the natural map $R\rightarrow R_P$.
Since $M$ is an $R_P$-module, $0_P\subseteq\mathrm{ann}_R(M)$. It is obvious that $0_P\subseteq A$. 
On the other hand, $AR_P$ is the nilradical of $R_P$. It follows that $\mathrm{rad}(0_P)=A$. Hence we get
that $\mathrm{Supp}\ M\subseteq V(\mathrm{ann}_R(M))\subseteq V(0_P)= V(A)$. If $B$ is another minimal prime, it is
obvious that $V(A)\cap V(B)=\emptyset$.

$(1)\Rightarrow (2)$. If $P$ is a prime ideal then $R_P=\mathrm{End}_R(R_P)$. It
follows that $P$ contains only one minimal prime ideal. \qed 
\begin{proposition} \label{P:Bind}
Let $R$ be a ring. The following conditions are
  equivalent:
\begin{enumerate}
\item For each $R$-algebra $S$ and for each indecomposable left $S$-module $M$, there is only one minimal prime ideal $A$ such that
  $\mathrm{Supp}\ M\subseteq V(A)$.
\item Each prime ideal contains a unique minimal prime ideal and $\mathrm{QMin}\ R$ is totally disconnected.
\end{enumerate}
\end{proposition}
\textbf{Proof.} $(1)\Rightarrow (2)$. By proposition~\ref{P:Bsupp} each
prime ideal contains a unique minimal prime ideal. Let $P\in\mathrm{QMin}\ R$ and $C$ its connected component.
There exists an ideal $A$ such that $V(A)=\lambda^{\leftarrow}(C)$.
Then $V(A)$ is  connected. It follows that
$R/A$ is indecomposable. So $V(A)=V(P)$ and $C=\{P\}$.

$(2)\Rightarrow (1)$. Let $S$ be an $R$-algebra and $M$ be an
indecomposable left $S$-module. Let $P\in\mathrm{Supp}\ M$, $A=\lambda(P)$,
$P'\in\mathrm{Spec}\ R\setminus V(A)$ and $A'=\lambda(P')$. Since $\mathrm{QMin}\ R$ is totally disconnected,
there exists an idempotent $e\in A\setminus A'$. We easily deduce that $e\in P\setminus P'$. Now we do as in the 
proof of Proposition~\ref{P:Gind} to conclude.
\qed

\end{document}